\documentclass[12pt]{amsart}

\title{On the geometry of a Picard modular group}

\author{Martin Deraux}

\date{May 23, 2022}

\newcommand{\F}{\mathbb{F}}
\newcommand{\Z}{\mathbb{Z}}
\newcommand{\R}{\mathbb{R}}
\newcommand{\Q}{\mathbb{Q}}
\newcommand{\C}{\mathbb{C}}
\newcommand{\cO}{\mathcal{O}}
\newcommand{\K}{\mathbb{K}}

\usepackage{graphicx}
\usepackage{a4wide}
\usepackage{url}

\newtheorem{prop}{Proposition}[section]
\newtheorem{thm}{Theorem}[section]
\def\cqfd{\mbox{}\nolinebreak\hfill$\Box$\medbreak\par}
\newenvironment{pf}{\noindent\textbf{Proof:}}{\cqfd}
\theoremstyle{definition}
\newtheorem{dfn}{Definition}[section]
\newtheorem{rk}{Remark}[section]
\newtheorem{algo}{Method}[section]

\usepackage{fancyvrb}
\usepackage{fvextra}

\begin{document}
\maketitle

\begin{abstract}
  We study geometric properties of the action of the Picard modular
  group $\Gamma=PU(2,1,\cO_7)$ on the complex hyperbolic plane
  $H^2_\C$, where $\cO_7$ denotes the ring of algebraic
  integers in $\Q(i\sqrt{7})$. We list conjugacy classes of maximal
  finite subgroups in $\Gamma$ and give an explicit description of the
  Fuchsian subgroups that occur as stabilizers of mirrors of complex
  reflections in $\Gamma$. As an application, we describe an explicit
  torsion-free subgroup of index $336$ in $\Gamma$.
\end{abstract}

\section{Introduction}

The goal of this paper is to study some detailed geometric features of
the Picard modular group $PU(J,\cO_7)$, where $\cO_7$ denotes the ring
of algebraic integers in $\Q(i\sqrt{7})$, and
\begin{equation}\label{eq:form}
  J=\left(\begin{matrix}
      0 & 0 & 1\\
      0 & 1 & 0\\
      1 & 0 & 0
\end{matrix}\right).
\end{equation}
Since $J$ has signature $(2,1)$, $PU(J)$ is isomorphic to the
group $PU(2,1)$ of holomorphic isometries of the complex hyperbolic
plane $H^2_\C$ (which is a symmetric space with $\frac{1}{4}$-pinched
sectional curvature).

It is a standard fact that $U(J,\cO_7)=U(J)\cap GL(3,\cO_7)$ is a
lattice in $U(J)$, i.e. a discrete subgroup such that the quotient of
$H^2_\C$ under the action of $U(J,\cO_7)$ has finite volume for the
symmetric Riemannian metric. From this point on, to simplify notation,
we write
$$
X=H^2_\C,\quad G=PU(J),\quad \Gamma=PU(J,\cO_7).
$$

The first explicit information about $\Gamma$ was an explicit finite
generating set, see Zhao~\cite{zhao}. A standard way to find such a
generating set is to work out a fundamental domain for the action,
because the side-pairing maps of a fundamental domain generate the
group. In principle, one could use the standard construction of Ford
domains, but the Ford domain for $\Gamma$ turns out to be very
complicated. In fact Zhao used coarser information than the actual
Ford domain.

More recently, his method was refined by Mark and
Paupert~\cite{mark-paupert} to get an actual presentation (generating
set plus defining relations for the group), still without working out
an explicit fundamental domain. The rough idea is to use a coarse
fundamental domain, i.e. a set $\Omega$ such that $\Gamma \Omega=X$,
and such that
\begin{equation}\label{eq:T}
  T=\{\gamma\in\Gamma:\Omega \cap \gamma\Omega\neq\emptyset\}
\end{equation}
is finite. In that case, it is easy to see that $T$ generates
$\Gamma$, and there is a simple way to write defining relations for
the group (see section~13.4 of~\cite{humphreys} for instance). In the
discussion that follows, we assume we have a coarse fundamental domain
and an explicit finite set $T$ as above.

The goal of the present paper is to push the Mark-Paupert techniques a
bit further, and give a list of the conjugacy classes of torsion
elements (and more generally of maximal finite subgroups) in
$\Gamma$. See section~\ref{sec:results} for the results.

This gives us some detailed information about the local structure of
the quotient orbifold $X/\Gamma$; as far as I know, its global
structure is still not understood. Note that, contrary to what was
incorrectly stated in~\cite{deraux-klein}, the group $\Gamma$ is
\emph{not} the same as the sporadic group
$\mathcal{S}(\bar\sigma_4,\infty)$, which is contained in $PU(H,\cO_7)$
for another Hermitian form. In fact the group
$\mathcal{S}(\bar\sigma_4,\infty)$ has trivial abelianization, whereas
$\Gamma$ has abelianization $\Z/2\Z$. This was our initial
motivation for studying the group $\Gamma$ in detail.

We will also determine the mirror stabililizers for the two conjugacy
classes of complex reflections in $\Gamma$ (see
section~\ref{sec:mirror-stabs}; it is a standard fact that these are
Fuchsian subgroups, but it is not exactly obvious how to describe
these stabilizers explicitly (generating set, signature). We find that
one stabilizer has a single cusp, whereas the other has two (the
latter case makes for much more complicated computations).

From our list of torsion conjugacy classes, we deduce the existence of
a torsion-free subgroup of index 336 in $\Gamma$, which is a principal
congruence subgroup (the kernel of the reduction modulo the ideal
$\langle i\sqrt{7}\rangle$). This is done in
section~\ref{sec:torsion-free-subgroup}.

Note that the methods used in this paper work for all 1-cusped Picard
groups (this happens for slightly more values of $d$, namely
$d=1,2,3,7,11,19,43,67,163$), but they require much heavier
computation. Prior to this work, presentations were worked out in the
literature only for $d=1,2,3,7,11$,
see~\cite{falbel-parker},~\cite{falbel-francsics-parker},
~\cite{mark-paupert},~\cite{polletta},
~\cite{ghoshouni-heydarpour-2},~\cite{ghoshouni-heydarpour-11}. Other
values of $d$ are treated in~\cite{deraux-xu}.

\noindent
\textbf{Acknowledgements:} The author wishes to thank Matthew Stover
for useful discussions related to this paper, as well as the anonymous
referee for his/her
careful reading of the first version of the manuscript.\\

\section{Complex hyperbolic geometry and Ford domains} \label{sec:cxhyp}

In this section we give a brief sketch of the geometry of the complex
hyperbolic plane, mainly to set up notation (we follow the notation
in~\cite{mark-paupert} quite closely). For much more detail, the
standard reference is~\cite{goldman-book}.

We work in homogeneous coordinates $v=(v_1,v_2,v_3)$ and write
$\langle v,w\rangle = w^* J v$, $||v||^2=\langle v,v\rangle$. As a
set, the complex hyperbolic plane $H^2_{\mathbb{C}}$ is the set of
complex lines in $\C^3$ that are spanned by a negative vector (i.e. a
vector $v\in\C^3$ with $||v||^2<0$). This set is contained in the
affine chart $v_3\neq 0$ of $P^2_{\C}$, where $v$ can be represented
as $(\frac{v_1}{v_3},\frac{v_2}{v_3},1)=(z_1,z_2,1)$, with
$2\Re(z_1)+|z_2|^2<0$. The distance function in complex hyperbolic
space is given by a simple formula in homogeneous coordinates, namely
$$
\cosh\left(\frac{1}{2}d([v],[w])\right)=\frac{ |\langle v,w \rangle| }{\sqrt{\langle v,v \rangle\langle w,w \rangle}},
$$
where $[v]=\C v$ denotes the complex line spanned by $v$. The boundary
at infinity of the complex hyperbolic plane, which consists of complex
lines spanned by null vectors (i.e. vectors $v\in\C^3$ with
$||v||^2=0$), is almost entirely contained in the affine chart
$z_3\neq 0$, only one point is missing, namely $q_\infty=(1,0,0)$. We
usually refer to that point as the point at infinity.

Rather than the affine coordinates $(z_1,z_2)$ described above, it is
convenient to use horospherical coordinates $(z,t,u)$, $z\in\C$,
$t,u\in\R$, $u>0$, defined by
$$
2z_1+|z_2|^2=it-u.
$$

Using these coordinates, the hypersurfaces defined by taking $u$ to be
a fixed positive constant are horospheres based at the point at
infinity. Points with $u=0$ give the boundary at infinity
$\partial_\infty H^2_\C$ of $H^2_\C$ (minus the point at infinity,
corresponding to $q_\infty=(1,0,0)$).

Each of these points has a unique representative of the form
$$
\left(
  \begin{matrix}
    \frac{-|z|^2+it}{2}\\
    z\\
    1
  \end{matrix}\right),
  $$
with $z\in\C$ and $t\in\R$.

Given a group $\Gamma\subset PU(J)$, we denote by $\Gamma_\infty$ the
stabilizer of $q_\infty=(1,0,0)$ in $\Gamma$, i.e. the subgroup of
matrices that have $q_\infty$ as an eigenvector. We start by
describing the full stabilizer in $PU(J)$.

The unipotent stabilizer of $q_\infty$ in $U(J)$ consists of the
matrices of the form
$$
T(z,t)=\left(\begin{matrix}
  1 & -\bar z & \frac{-|z|^2+it}{2}\\
  0 &    1    &     z\\
   0 &    0    &     1
\end{matrix}
\right),
$$
and it acts simply transitively on $H^2_\C\setminus\{\infty\}$. This
gives $H^2_\C\setminus\{\infty\}$ the structure of a group, usually
called the Heisenberg group. In terms of the coordinates
$(z,t)\in\C\times\R$, the Heisenberg group law is the following:
\begin{equation}
  \label{eq:heis-law}
  (z,t)\star(z',t')=(z+z',t+t'+2Im(z\bar z')).
\end{equation}

Non-unipotent parabolic elements are
usually called twist-parabolic elements; they can be written as
$RT(z,t)$ where $R=diag(1,\zeta,1)$, $|\zeta|=1$, i.e.
\begin{equation} \label{eq:twist-parabolic}
RT(z,t)=\left(\begin{matrix}
  1  & -\bar z  & \frac{-|z|^2+it}{2}\\
  0  &    \zeta &     z\\
   0 &    0     &     1
\end{matrix}
\right).
\end{equation}

Recall that $q_\infty=(1,0,0)$.
\begin{dfn}\label{dfn:ford}
  Let $\Gamma$ be a discrete subgroup of $G$. The Ford domain for
  $\Gamma$ is defined as
  $$
  F_\Gamma = \{ [x]\in H^2_\C\ :\ |\langle x,q_\infty\rangle| \leq |\langle x,\gamma q_\infty\rangle| \textrm{ for all }\gamma\in\Gamma \}
  $$  
\end{dfn}
In this definition, $\gamma q_\infty$ stands for
$\tilde{\gamma} q_\infty$ for any lift $\tilde{\gamma}\in U(J)$ of
$\gamma\in \Gamma\subset PU(J)$. Since lifts of a given element differ
by multiplication by a complex number with modulus one, the
inequalities in Definition~\ref{dfn:ford} are independent of the lift
chosen.

It is easy to see that $F_\Gamma$ is invariant under the action of
$\Gamma_\infty$; indeed, if $\alpha\in\Gamma_\infty$ and
$x\in F_\Gamma$, then
$$
|\langle \alpha x,q_\infty\rangle| = |\langle
x,\alpha^{-1}q_\infty\rangle| = |\langle x,q_\infty\rangle| \leq
|\langle x,\alpha^{-1}\gamma q_\infty\rangle| = |\langle \alpha
x,\gamma q_\infty\rangle|,
$$
so $\alpha(x)\in F_\Gamma$.
In particular, $F_\Gamma$ cannot always be a fundamental domain for
the action of $\Gamma$. It is a standard fact that it is a fundamental
domain if (and only if) $\Gamma_\infty$ is trivial (see section~9.5
of~\cite{beardon} for a proof in the complex 1-dimensional case).

When $\Gamma_\infty$ is not trivial, in order to get a fundamental
domain for $\Gamma$, we select a fundamental domain $P$ for the action
of $\Gamma_\infty$ in the Heisenberg group, and consider the cone
$C_P$ with base $P$ in horospherical coordinates:
$$
C_P = \{(z,t,u)\ :\ (z,t)\in P, u>0 \}.
$$
Then we have:
\begin{prop}
  $F_\Gamma\cap C_P$ is a fundamental domain for $\Gamma$.
\end{prop}

Note that even though $F_\Gamma\cap C_P$ need not have well-defined side-pairing
maps, the sides of the Ford domain $F_\Gamma$ are paired. In
fact, given $\gamma\in\Gamma\setminus\Gamma_\infty$, write
$$
I(\gamma) = \{ [x]\in H^2_\C\ :\ |\langle x,q_\infty\rangle| \leq |\langle x,\gamma q_\infty\rangle| \}.
$$
It follows easily from the definition that
$\gamma^{-1}(I(\gamma))=I(\gamma^{-1})$, hence the set
$F_\Gamma\cap I(\gamma)$, if it is a side (i.e. if it has dimension
3), must be paired with $F_\Gamma\cap I(\gamma^{-1})$ by
$\gamma^{-1}$.

The set $I(\gamma)$, $\gamma\in \Gamma\setminus\Gamma_\infty$ can be interpreted as a
bisector (locus equidistant of two points in $H^2_\C$), or as an
isometric sphere (locus of points where the Jacobian of the
transformation is 1), or as a metric sphere for the so-called extended
Cygan distance, defined by
\begin{equation}\label{eq:cygan-distance}
d_{\textrm{Cygan}}((z,t,u),(z',t',u'))=\left( (|z-z'|^2 + |u-u'|)^2 + |t-t'+2\textrm{Im}(z\bar z')|^2 \right)^{1/4}.  
\end{equation}
Using this distance, we can describe $I(\gamma)$ as the Cygan sphere
with center $\gamma(\infty)$, and radius $\sqrt{2/|a_{31}|}$, where
$A=(a_{jk})_{j,k=1,2,3}$ is a matrix representative of $\gamma$
(see~\cite{mark-paupert} for instance).

It can be useful also to have an explicit expression for the Cygan
sphere of radius $r$ centered at the point $(c,d)\in\C\times\R$ in
the Heisenberg group, namely it consists of points with horospherical
coordinates $(z,t,u)\in\C\times\R\times\R_+$ satisfying
\begin{equation}\label{eq:cygan-sphere}
 ( |z-c|^2 + u )^2 + |t-d + 2\textrm{Im}(z\bar c)|^2  = r^4.
\end{equation}
The basic observation is that if $(z,t,u)$ is in that sphere, then
$|z-c|\leq r$ so the $\C$-component must be contained in the Euclidean
disk of radius $r$ centered at $c$. We also have the basic estimates
$u\leq r^2$, and a slightly less efficient estimate for $t$ given by 
$$
  |t-d + 2\textrm{Im}(z\bar c)| \leq r^2,
$$
which gives an estimate for the range of values of $t$ for any fixed
value of $z$.

We now focus on the special case of $\Gamma=PU(J,\cO_7)$, and review
some results of~\cite{mark-paupert} giving an explicit description of
$\Gamma_\infty$ (see also~\cite{paupert-will}).

We define
\begin{equation}\label{eq:tsl}
  T_1=T(1,\sqrt{7}),\quad T_\tau=T(\tau,0),\quad T_v=T(0,2\sqrt{7}),
\end{equation}
where $\tau=\frac{1+i\sqrt{7}}{2}$.

Since we consider $\Gamma=PU(J,\cO_7)$, in the definition of
twist-parabolic elements given in equation~\eqref{eq:twist-parabolic},
we only allow $\zeta$ to be a unit in $\cO_7$, i.e. $\pm 1$. Here and
in what follows, we define
\begin{equation}\label{eq:R}
  R=\textrm{diag}(1,-1,1)
\end{equation}

\begin{figure}[htbp]
  \centering
  \includegraphics[width=0.3\textwidth]{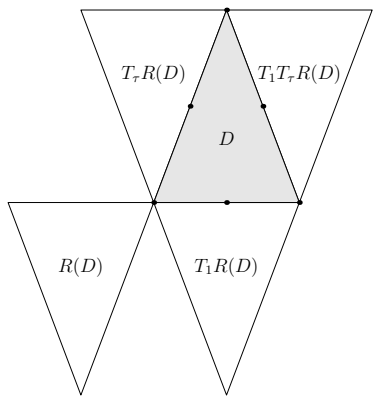}
  \caption{Action of the cusp group on the $\C$-factor of the Heisenberg group $\C\times\R$.}\label{fig:triangle}
\end{figure}

\begin{figure}[htbp]
  \centering
  \includegraphics[height=0.2\textheight]{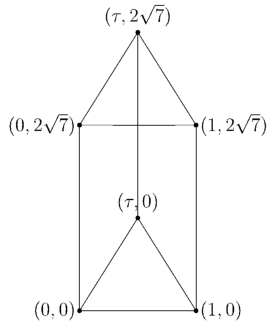}
  \caption{$P=D\times[0,2\sqrt{7}]$ is a fundamental domain for the action of $\Gamma_\infty$ on $\C\times\R$}\label{fig:prism}
\end{figure}

In the $\C$-factor of the Heisenberg group $\C\times\R$, $R$ acts as
$z\mapsto -z$, and $T_1,T_\tau$ act as translations by $1$ and $\tau$
respectively (whereas $T_v$ act trivially). Hence the triangle $D$
which is the Euclidean convex hull of $0,1$ and $\tau$ gives a
fundamental domain for the action on $\C$. Note that $T_1R$,
$T_\tau R$ and $T_1T_\tau R$ act in the $\C$ factor as half-turns
fixing the midpoints of the sides of the triangle, see
Figure~\ref{fig:triangle}.
The group $\Gamma_\infty$ is actually generated by $T_1,T_\tau$ and
$R$, see~\cite{paupert-will} for instance (note that
$T_v=[T_\tau,T_1]$). Since $T_v$ acts as a translation by $2\sqrt{7}$
in the $t$-coordinate, it should be quite clear that the prism
$P=D\times[0,2\sqrt{7}]$ is a fundamental domain for the action of
$\Gamma_\infty$ in $\C\times\R$, see Figure~\ref{fig:prism}.
Note that $T_\tau R$ and $T_1T_\tau R$ are side-pairing maps for $P$
(in fact these are complex reflections), and $T_v$ as well (it gives
the vertical translation pairing the top and bottom triangles of the
prism). However $T_1 R$ is not a side pairing-map (it is given by a
glide-reflection), but this will be inconsequential in the present
paper.

The domain $P$ is chosen to have affine sides in Heisenberg
coordinates (since the Heisenberg group acts on itself by affine
transformations, see formula~\eqref{eq:heis-law}). It can also be
adjusted to have well-defined side-pairing maps (this requires
subdividing its sides into smaller polygons), but we will not do this,
since we do not need it for any of the methods used in this paper.

A few Ford domains in $H^2_\C$ have been studied explicitly,
see~\cite{deraux-family},~\cite{parker-will} among others. Such
domains are virtually present in~\cite{zhao} and~\cite{mark-paupert},
but it is very complicated to determine their combinatorial structure
in detail (or even to describe it on paper!)

\begin{figure}[htbp]
  \centering
  \includegraphics[width=0.6\textwidth]{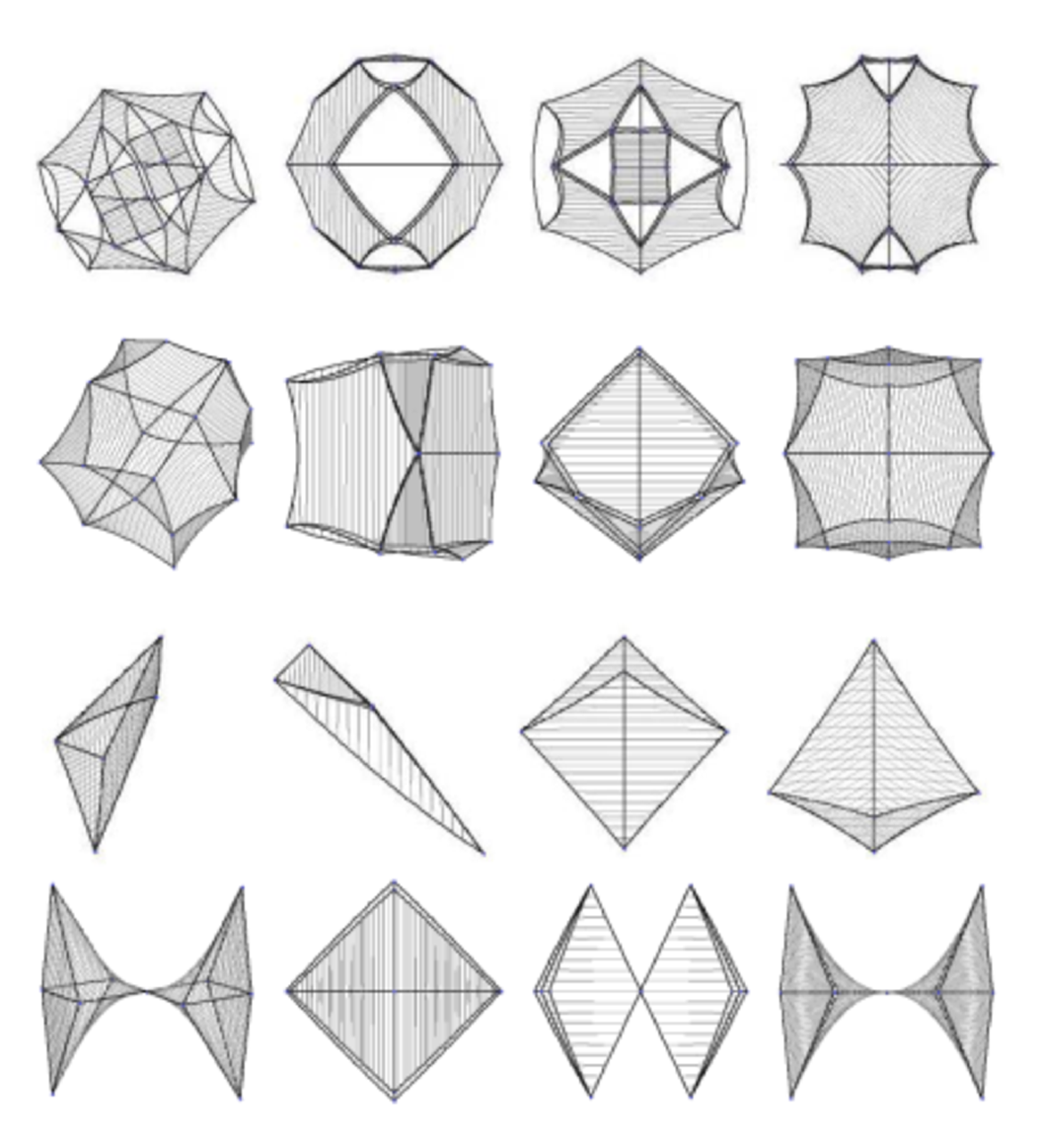}
  \caption{Sides of the Ford domain for $PU(J,\cO_7)$ (four views for each side).}
  \label{fig:ford-O7}
\end{figure}

To give a rough idea, pictures of (representatives of the isometry
type of) the sides of the Ford domain for $\Gamma=U(J,\cO_7)$ are given
in Figure~\ref{fig:ford-O7}. Needless to say, these pictures will not
be used anywhere in the paper, but they should give an idea of the
intricacy of the combinatorics.

\section{Virtual fundamental domain, algorithms} \label{sec:algos}

Recall that a coarse fundamental domain for a discrete group $\Gamma$
is a set $\Omega$ such that $\Gamma \Omega=X$, and such that
$T=\{\gamma\in\Gamma:\Omega \cap \gamma\Omega\neq\emptyset\}$ is
finite. In that case, one can prove that $T$ generates $\Gamma$, and
write a explicit group presentation in terms of these generators (see
section 13.4 in~\cite{humphreys}).

Of course a fundamental domain is a coarse fundamental domain, and
this is the kind of coarse fundamental domain we will (virtually) use
in this paper. Specifically, we will use $\Omega=F_\Gamma\cap C_P$,
where $C_P$ is the cone to the point at infinity in horospherical
coordinates, and $P$ is the prism described in~\cite{mark-paupert},
i.e.
$$
\Omega = \left\{ (a+b\tau,t)\ : 0\leq a,b,a+b\leq 1,\ 0\leq t\leq 2 \right\}.
$$
By making virtual use of this domain, we mean that, rather than
studying the combinatorics of $\Omega$, we will only use
algorithmic procedures that allow us to do the following.
\begin{algo} \label{algo:key}
\begin{enumerate}
\item Determine whether a given algebraic point $[x]\in H^2_\C$ is
  inside $\Omega$, and list the sides of $\Omega$ that contain it;
\item If $[x]$ is not in $\Omega$, find an element $\gamma\in\Gamma$
  such that $\gamma[x]\in\Omega$.
\end{enumerate}  
\end{algo}
By an algebraic point $[x]\in H^2_\C$, we mean a point that can be
represented in homogeneous coordinates by a vector $x\in\C^3$ with
algebraic coordinates. Most algebraic points we consider are actually
$\K$-rational (i.e they are represented by a vector in $\K^3$), but
not all; indeed, the isolated fixed points of a regular elliptic
isometry is an eigenvector of a matrix with entries in $\K$, but their
eigenvalue need not be in $\K$ in general.

We briefly explain why such methods can be implemented effectively. The
main difficulty is that $[x]\in \Omega$ is defined in terms of an
infinite set of inequalities.

We first make use of some results in~\cite{mark-paupert}, which allow
us to reduce the list of $\gamma\in\Gamma$ that occur
definition~\eqref{dfn:ford}. This reduction is based on arithmetic
considerations.

Here an in what follows, we write $\K=\Q(i\sqrt{7})$, and $\cO_7$ for
the ring of algebraic integers in $\K$.
\begin{dfn}
  A $\K$-rational point in $P^2_\C$ is a point that can be represented
  in homogeneous coordinates as a vector with entries in $\cO_7$.
\end{dfn}
A given a rational point has an essentially unique primitive
representative in $\cO_7^3$, i.e. a vector whose entries have greatest
common divisor equal to 1 (see Lemma~1 of~\cite{mark-paupert} or
Lemma~3 of~\cite{polletta}). The only freedom we have to change such a
vector is to scale it by unit, and the only units in $\cO_7$ are
$\pm 1$.

By the square norm of a point in $P^2 \C$, we mean $||v||^2=v^*Jv$ for
any primitive representative $v\in\cO_7^3$. Points (resp. vectors) with
square norm 0 are called null points (resp. vectors). We will use the
following result (see Lemma~2 in~\cite{mark-paupert}).
\begin{prop}\label{prop:1cusp}
  For every $A\in U(J,\cO_7)$, the first column of $A$ (which
  represents $A(\infty)$) is a primitive null vector in
  $\cO_7^3$. Conversely, given any primitive null vector $v\in \cO_7^3$,
  there exists an $A\in U(J,\cO_7)$ whose first column is $v$.
\end{prop}
The converse direction says that the Picard group $\Gamma=U(J,\cO_7)$
acts transitively on $\K$-rational null vectors; this is can be seen
as a reformulation of the statement that the quotient $X/\Gamma$ has
exactly one cusp, which is a consequence of the fact that $\K$ has
class number 1 (see~\cite{zink} or section~4.1
in~\cite{stover-cusps}). The theoretical result in
Proposition~\ref{prop:1cusp} is unfortunately difficult to make
effective, as discussed in~\cite{mark-paupert}.

We review the following definition (see~\cite{mark-paupert}).
\begin{dfn}
  \begin{enumerate}
  \item The \emph{depth of a $\K$-rational point}
    $[v]\in \partial_\infty H^2\C\setminus\{\infty\}$ is given by
    $|v_3|^2$, where $v=(v_1,v_2,v_3)$ is any primitive integral
    representative of $[v]$.
  \item The \emph{depth of a matrix}
    $A\in\Gamma\setminus\Gamma_\infty$ is the depth of $A(\infty)$.
\end{enumerate}  
\end{dfn}
The notion of depth is important, in part because of the following
Proposition (see Proposition~4.3 of~\cite{kim-parker}).
\begin{prop}
  Let $A\in\Gamma$. Then the Cygan center of the isometric sphere
  $I(A)$ is represented by $A(\infty)$, and its radius is given by
  $\sqrt{2/d}$, where $d$ is the depth of $A$.
\end{prop}

Note that the depth of any $A\in \Gamma$ is a rational integer, but
not every integer occurs as the depth of an element $A\in\Gamma$ (only
integers that are norms in $\cO_7$ occur). For example, the first five
numbers that occur as depths of elements in the Picard group
$\Gamma$ are $1,2,4,7,11$.

If there exists an element of a given depth, then there exist
infinitely many, since pre- or post-composition with any element of
$\Gamma_\infty$ does not change the height. However, in a fixed
bounded region of the Heisenberg group
$\C\times\R\simeq \partial_\infty H^2_\C\setminus\{\infty\}$, there
are only finitely many points of a given depth. This is true in
particular in our fundamental domain $P$ for the action of
$\Gamma_\infty$. A list of $\K$-rational points of small depth in $P$
is given in ~\cite{mark-paupert}.

The following result is a consequence of the covering depth estimate
given in~\cite{mark-paupert}.
\begin{thm}
  The isometric spheres of elements of depth $\geq 11$ do not
  intersect the Ford domain $F_\Gamma$.
\end{thm}
\begin{rk}
  It turns out that the isometric spheres of the elements of depth 7
  \emph{do} intersect the Ford domain, but in lower-dimensional facets
  only, i.e. they are not needed in the
  definition~\eqref{dfn:ford}. This fact is quite painful to prove
  however, and we will not use it in the sequel.
\end{rk}

We will also use the following list of 14 elements in the group, which
is a slight modification of the list given in~\cite{mark-paupert} (we
want the set of side pairing maps to be closed under inversion of
matrices).
\begin{equation*}
  \begin{array}{c}
  A_1 =
  \left(
    \begin{matrix}
      0 & 0 & 1\\
      0 & -1 & 0\\
      1 & 0 & 0
    \end{matrix}
  \right),\quad
  A_2 =
  \left(
    \begin{matrix}
      2 & -\tau & 1-3\tau\\
      \bar\tau & 0 & -2-\tau\\
      -\tau & -1 & -3+\tau
    \end{matrix}
  \right),\quad
  A_3 = A_2^{-1},\quad A_4=A_2^{-2},\quad A_5 = A_4^{-1}
  \\
    A_6 =
  \left(
    \begin{matrix}
      i\sqrt{7} & 0 & 4\\
      0 & 1 & 0\\
      2 & 0 & -i\sqrt{7}
    \end{matrix}
  \right),\quad
  A_7 =
  \left(
    \begin{matrix}
      -\bar\tau & 1 & 1\\
      \tau & 0 & 1\\
      2 & \bar\tau & -\tau
    \end{matrix}
  \right),\quad
  A_8 =
  \left(
    \begin{matrix}
      1 & 2-\tau & -2\\
      -1-\tau & -3 & 1+\tau\\
      -2 & -2+\tau & 1
    \end{matrix}
  \right),\quad
  \\
    A_9 =
    \left(
    \begin{matrix}
      -1 & 0 & i\sqrt{7}\\
      0 & 1 & 0\\
      i\sqrt{7} & 0 & 6
    \end{matrix}
  \right),\quad
  A_{10} =
  \left(
    \begin{matrix}
      -2 & 0 & i\sqrt{7}\\
      0 & 1 & 0\\
      i\sqrt{7} & 0 & 3
    \end{matrix}
  \right),\quad
  A_{11}=A_{10}^{-1},\quad
    \\
  A_{12} =
    \left(
    \begin{matrix}
      -4 & 0 & 3i\sqrt{7}\\
      0 & 1 & 0\\
      i\sqrt{7} & 0 & 5
    \end{matrix}
  \right),\quad
  A_{13}=A_{12}^{-1},\quad
  A_{14}=A_{9}^{-1}
  \end{array}  
\end{equation*}

The basic fact we will use from~\cite{mark-paupert} is the following.
\begin{thm} \label{thm:bound-height}
  Let $\gamma\in \Gamma\setminus\Gamma_\infty$ be such that $I(\gamma)$
  intersects the Ford domain. Then there exists
  $\alpha,\beta\in\Gamma_\infty$ and $j\in\{1,\dots,14\}$ such that
  $\gamma=\alpha A_j \beta^{-1}$.
\end{thm}
The remaining difficulty is that $\Gamma_\infty$ is an infinite group,
and we cannot check infinitely many inequalities with a computer; we
now explain how to handle this difficulty.

In sections~\ref{sec:basic-algo} through~\ref{sec:isotropy}, we sketch
the general algorithms we will use. We then list the specific results
(conjugacy classes of torsion elements, maximal finite subgroups) for
the group $U(J,\cO_7)$ in section~\ref{sec:results}.

\subsection{Basic algorithms} \label{sec:basic-algo}

Recall that we are after algorithmic methods stated
in Method~\ref{algo:key}. Let us assume $x\in\C^3$ has algebraic coordinates,
and $||x||^2<0$. We can compute the (algebraic) horospherical
coordinates $(z,t,u)$ of $x$, using
$$
z = \frac{x_2}{x_3}, it-u = 2\frac{x_1}{x_3}+|\frac{x_2}{x_3}|^2.
$$
We would like to find a cusp element $\gamma\in\Gamma_\infty$ such
that $\gamma(x)$ has horospherical coordinates $(z',t',u')$ with
$(z',t')\in P$. Recall that $P\subset \C\times\R$ is
$T\times[0,2\sqrt{7}]$ where $T$ is the convex hull of $0,1$ and
$\tau$.

First solve
$$
a+b\tau = z\leftrightarrow \left\{\begin{array}{l} a+\frac{b}{2} = \textrm{Re}(z)\\ b\sqrt{7}=\textrm{Im}(z)\end{array}\right.
$$
for real algebraic $a,b$, and compute the floors
$k = \lfloor a\rfloor$ and $l=\lfloor b\rfloor$,
$\epsilon=\lfloor a+b\rfloor$. By applying $T_1^{-k}T_\tau^{-l}$ for
suitably chosen $k,l\in\Z$, we may assume the horospherical
coordinates $(z,t,u)$ satisfy $z\in T$. Computing
$\lfloor t/(2\sqrt{7})\rfloor$, we get a power of $T_v$ that brings
the $t$-coordinate in $[0,2\sqrt{7}]$.

In other words, we may assume that $(z,t)\in P$. Now we use the
following.
\begin{prop}\label{prop:cusp-shift}
  Let $j\in \{1,\dots,14\}$. There is an explicit set
  $E_j\subset \Gamma_\infty$ such that for all
  $\gamma\in\Gamma_\infty\setminus E_j$, we have
  $\gamma(I(A_j))\cap C_P=\emptyset$.
\end{prop}
Without the ``explicit'' request, this is a consequence of
discreteness of $\Gamma_\infty$. The explicit character follows from
elementary estimates using the values of the Cygan radius of $I(A_j)$
(for details on these estimates, see~\cite{deraux-xu}).

Using Proposition~\ref{prop:cusp-shift}, we can list all the isometric
spheres $\alpha(I(A_j))$ that contain $[x]$ (or such that $[x]$ is in
the interior of the corresponding Cygan ball). This gives method~(1).

For method~(2), we first proceed as above to bring the Heisenberg
coordinate to $P$, then find the Cygan spheres $\alpha(I(A_j))$ that
contain $x$, and apply the isometry $(\alpha A_j)^{-1}$ to $x$. This
decreases the number of Cygan balls defining the Ford domain that
contain $[x]$ in their interior. Then repeat the preceding procedure
until that number is 0.

\subsection{Conjugacy classes of torsion elements} \label{sec:algo-torsion}

Clearly every torsion element is conjugate to a torsion element that
fixes a point in the fundamental domain $\Omega = F_\Gamma\cap
C_P$. Moreover, the dicreteness of $\Gamma$ implies that there exists
a precisely invariant horoball, i.e. a horoball $B$ based at
$q_\infty$ such that is invariant under $\Gamma_\infty$, and
$\gamma(B)\cap B=\emptyset$ for every
$\gamma\in\Gamma\setminus\Gamma_\infty$. Moreove, the fact that
$\Gamma$ is a lattice implies that $\Gamma_\infty$ acts cocompactly on
every horosphere based at $q_\infty$.

It follows that there are actually finitely many torsion elements
fixing a point in $\Omega$; we now explain a method for listing these
elements.

Suppose $\gamma\in\Gamma$ has finite order. There are two
possibilities, either $\gamma\in\Gamma_\infty$ or not.

\subsubsection{$\gamma\in\Gamma_\infty$}

In this case, $\gamma$ must be a complex reflection, and we may assume
that its mirror meets the Heisenberg group $\C\times\R$ in a vertical
line that intersects the boundary of the prism $P$. Projecting to the
$\C$-factor, it is easy to see that the boundary at infinity of the
mirror is given by $\{(z_0,t):t\in\R\}$ for $z_0=0$, $1/2$, $\tau/2$,
$(1+\tau)/2$. This implies that $\gamma$ has one of the following
forms:
$$
RT_v^k,\quad T_1RT_v^k,\quad T_\tau RT_v^k,\quad T_1T_\tau RT_v^k,
$$
for some $k\in\Z$ (we will give more details about this in
section~\ref{sec:cusp-torsion}). The only complex reflections of this
form are
$$
R,T_1R,T_1T_\tau R.
$$

\subsubsection{$\gamma\notin\Gamma_\infty$}

Suppose $\gamma$ fixes a point $[x]\in P$. It is easy to see from the
definition of isometric spheres that we must have
$[x]\in I(\gamma)\cap I(\gamma^{-1})$, and in particular
$$
I(\gamma)\cap I(\gamma^{-1})\neq\emptyset.
$$

Also, because the fixed point $[x]$ is in the Ford domain, $\gamma$
must have depth $\leq 7$, so $\gamma=\alpha A_j\beta^{-1}$ for some
$j\in\{1,\dots,14\}$ and $\alpha,\beta\in\Gamma_\infty$ (see
Theorem~\ref{thm:bound-height}).

Since we only want to list torsion elements up to conjugacy, we may
assume $\beta=Id$, i.e. $\gamma=\alpha A_j$, and we get
$$
\emptyset \neq I(\gamma)\cap I(\gamma^{-1}) = \alpha(I(A_j)) \cap I(A_j^{-1}).
$$
Note that we chose the set $\{A_1,\dots,A_{14}\}$ to be invariant
under the operation of taking inverses, so there is a
$k\in\{1,\dots,14\}$ such that $I(A_j^{-1})=I(A_k)$, and we must have
$$
\alpha(I(A_j))\cap A_k\neq \emptyset.
$$

\begin{prop} \label{prop:explicit-disjoint}
  There is an explicit finite set $T_{jk}$ such that for all
  $\alpha\in\Gamma_\infty\setminus T_{jk}$,
  $\alpha(I(A_j))\cap A_k= \emptyset$.
\end{prop}
The set can be made explicit by using elementary estimates using the
triangle inequality and the known radii of the Cygan spheres $I(A_j)$,
$I(A_k)$ (see equation~\eqref{eq:cygan-sphere}). For more details,
see~\cite{deraux-xu}.

From this, we get that every torsion element in
$\Gamma\setminus\Gamma_\infty$ is conjugate to an element of the
finite sets $T_{jk}$ for $j,k\in\{1,\dots,14\}$ (and we may restrict
to pairs $j,k$ such that $A_jA_k=Id$).

Moreover, given an element $\gamma\in T_{jk}$, there is an algorithm
to determine whether $\gamma$ has finite order. Indeed, the
eigenvalues of the matrix representative for $\gamma$ (which is unique
up to multiplication by $-Id$) are algebraic integers, and we can
determine whether they are roots of unity by examining their minimal
polynomial.

If the eigenvalues are all roots of unity, we can check whether or not
the matrix is diagonalizable by computing its minimal
polynomial. Hence we have an algorithm to do the following.
\begin{algo}
  Produce a finite set that contains a representative of every
  conjugacy class of torsion element in $\Gamma$.
\end{algo}

For elements with isolated fixed points, we can use methods of
section~\ref{sec:basic-algo}, we remove all elements whose fixed point
set is not in the fundamental domain $\Omega$, since they must be
conjugate to another element in the list.

\subsection{Eliminating redundancy} \label{sec:redundancy}

The list obtained by applying the method explained in
section~\ref{sec:algo-torsion} may have redundancies, in the sense
that some elements in the list may be conjugate to each other.

We now explain how to test whether two torsion elements $\gamma_1$,
$\gamma_2$ are conjugate to each other, assuming that they both have
an isolated fixed point. We call $x_j$ the isolated fixed point of
$\gamma_j$; by the methods of section~\ref{sec:basic-algo}, we
may assume $x_1,x_2$ are both in $\Omega$.

Suppose $\gamma_1,\gamma_2$ are conjugate in $\Gamma$, i.e there
exists $\alpha\in\Gamma$ such that
$\gamma_2=\alpha\gamma_1\alpha^{-1}$. Then $\alpha(x_1)=x_2$, and in
particular
$$
\Omega\cap\alpha(\Omega)\neq\emptyset.
$$
This implies that for some $j\in\{1,\dots,14\}$, $\alpha(I(A_j))$ must
intersect $P$, and there are finitely many choices of $\alpha$ such
that this is the case (as before, this can be made explicit).

Hence we only need to check finitely candidate conjugators $\alpha$ in
order to determine whether the two elements are conjugate, making this
special case of the conjugation problem solvable.

For pairs of complex reflections, it is more complicated to write a
general algorithm to test conjugacy (because the stabilizer of the
mirror of one such complex reflection is infinite).  In order to treat
the special case of $PU(J,\cO_7)$, it suffices to use the following two
observations
\begin{enumerate}
\item It two elements are conjugate, then we can find an element that
  conjugates them by listing all group elements, by listing words in a
  fixed generating set in increasing word length;
\item If two complex reflections have different order, or different
  Jordan forms, or different square norm for their (primitive) polar
  vector, then they are not conjugate.
\end{enumerate}
For a general group, we are likely to find pairs of elements with the
same rough conjugacy invariants (order, Jordan form, square norm of
polar vector) but the enumeration of the group fails finding a
conjugator (say because we run out of time or memory).

\subsection{Finding maximal finite subgroups} \label{sec:isotropy}

In this section, we explain how to determine maximal finite
subgroups. We assume that we have applied the methods of
section~\ref{sec:algos} successfully, and that we have a finite list
of all torsion elements whose fixed points contain a point in
$\Omega$.

The interesting maximal finite subgroups contain an element with an
isolated fixed point, since generic points on the mirror of a complex
reflection have a cyclic group as their stabilizer, generated by a
single complex reflection.

Now take the list of torsion elements with an isolated fixed point in
$\Omega$. For each such fixed point $[x]$, we use the methods of
section~\ref{sec:basic-algo} to find the list of Ford spheres and
sides of $C_P$ containing $[x]$.

Now build a graph whose vertex set is in bijection with the set of
these isolated fixed points, and join two vertices by a directed edge
if there is a side-pairing map sending one to the other (either coming
from the Ford domain, or from the prism $P$); note that the edges may
join a vertex to itself.

The conjugacy classes of maximal finite groups in $\Gamma$ are then
given by connected component of this graph (take the image of the
obvious representation of the fundamental group of the graph into
$\Gamma$).

For concreteness, we work out a couple of explicit examples.
\begin{enumerate}
\item 
Consider
$M=(RT_1IT_1^{-1})^2$ for instance (see the first entry in
Table~\ref{tab:order2}). Its isolated fixed point is given by $[v_1]$
where $v_1 = (-\bar\tau,0,1)$, which has horospherical coordinates
$[0,\sqrt{7},1]$. This point is on two sides of the cone $C_P$ (namely
the ones with side-pairing $T_\tau R$ and $T_1R$), as well as on three
Ford-Cygan spheres , namely
  $$
  I(A_6), T_1(I(A_1)),  T_1^{-1}T_v(I(A_1)).
  $$
  The side-pairing maps associated to these sides are given respectively by
  $$
  A_6, T_1A_1T_1^{-1}, T_1^{-1}T_vA_1T_v^{-1}T_1,
  $$
  and each of these elements actually fixes $[v]$.

  Using the side-pairings coming from $C_P$, we find two other points
  $[v_2]$, $[v_3]$ given $v_2=T\tau Rv_1$ and $v_3=T_1 Rv_1$. Each is
  on three Ford-Cygan spheres, which are simply the images of the above
  three Cygan spheres under $T_\tau R$ (or $T_1 R$).

  Hence the connected component of the graph containing $[v_1]$ is a
  triangle as in Figure~\ref{fig:graph1}. The stabilizer of $[v_1]$ is
  generated by $A_6, T_1A_1T_1^{-1}, T_1^{-1}T_vA_1T_v^{-1}T_1$
  together with the element
  $$
  (T_1R)^{-1}T_1T_\tau R T_\tau R = R.
  $$
  \begin{figure}[htbp]
    \centering
    \includegraphics[width=0.6\textwidth]{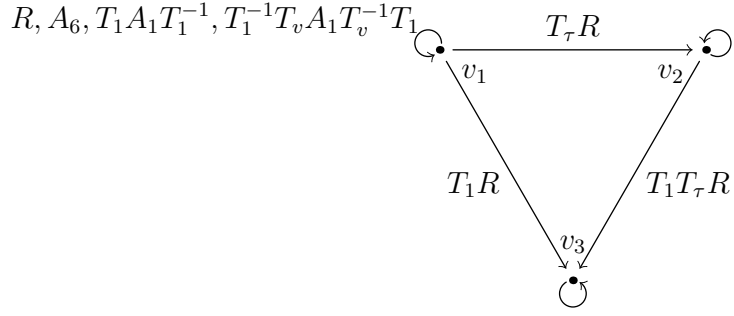}
    \caption{Cycle graph for the vertex $v_1$. We omit the Ford
      side-pairing maps that fix $v_2$ (resp. $v_3$), since these are
      simply conjugates of the elements in the label for the loop from
      $v_1$ to itself.}
    \label{fig:graph1}
  \end{figure}
  The linear group generated by
  $R, A_6, T_1A_1T_1^{-1}, T_1^{-1}T_vA_1T_v^{-1}T_1$ has order 16,
  and its subgroup of scalar matrices has order 2. In other words, the
  stabilizer of $[v_1]$ in $\Gamma$ has order 8.

  It is easy to check (for instance by enumerating the elements in the
  stabilizer) that the stabilizer contains four complex
  reflections.

  Let us denote $v_1,\dots,v_4$ polar vectors to the mirrors of these
  reflections, which we may and do chose to be primitive vectors in
  $\cO_7^3$. Explicit computation shows that, perhaps after permuting
  these vectors, we have
  $\langle v_1,v_1\rangle=\langle v_2,v_2\rangle=1$,
  $\langle v_1,v_1\rangle=\langle v_2,v_2\rangle=2$.
  
  This kind of information is gathered in columns 4 through 6 of
  Tables~\ref{tab:order2}-\ref{tab:order7} (see also
  Definition~\ref{def:jline}).
\item Consider the element
  $N=(T_\tau R)\cdot T_1R(T_1I)^2T_1^{-1}RT_1IT_1^{-1}\cdot
  (RT_\tau^{-1})$, which is the element of order 6 given in
  Table~\ref{tab:order6}, conjugated by $T_\tau R$ so that its
  isolated fixed point $[w]$ is in $\Omega$.

  The corresponding point is not on any side of the cone $C_P$, but it
  is on five Ford Cygan spheres, namely
  $$
  I(A_2), I(A_3), RT_\tau^{-1}(I(A_4)), RT_\tau^{-1}T_1^{-1}(I(A_5)), T_\tau(I(A_6)).
  $$  
  The point $[A_2^{-1}w]$ is inside the Ford domain, but it is not in
  the cone $C_P$, so we bring it back to the cone by a cusp
  element. It turns out $T_\tau R$ does the job, i.e.
  $T_\tau RA_2^{-1}w\in\Omega$. Note that the Cygan sphere for
  $(T_\tau RA_2^{-1})^{-1}=A_2\cdot RT_\tau^{-1}$ is of course the same as
  the one for $A_2$, since these two elements differ by
  pre-composition by a cusp element.

  Similar considerations show that (adjusted) side-pairing maps
  corresponding to the above five Cygan spheres are
  \begin{equation*}
    \begin{array}{c}
      T_\tau R \cdot A_2^{-1},\\
      T_1T_\tau R \cdot A_3^{-1},\\
      T_\tau^2 R \cdot (T_\tau A_6^{-1} T_\tau^{-1}),\\
      T_1 \cdot (R T_\tau^{-1} A_4^{-1} T_\tau R),\\
      T_1^{-1} \cdot (RT_\tau^{-1}T_1^{-1} A_5^{-1} T_1T_\tau R)\\
    \end{array}
  \end{equation*}
  and all these elements fix $[w]$. In other words, the relevant
  (connected component of the) cycle graph has a single vertex, and 5
  loops.

  The above five matrices generate a linear group of order 12, whose
  projectivization is a cyclic group of order 6. There is a complex
  reflection in this group, obtained by taking the third power of a
  generator (see Table~\ref{tab:order6}).
\end{enumerate}
\section{Results}\label{sec:results}

\subsection{Torsion in $\Gamma_\infty$} \label{sec:cusp-torsion}

\begin{figure}[htbp]
  \centering
  \includegraphics[height=6cm]{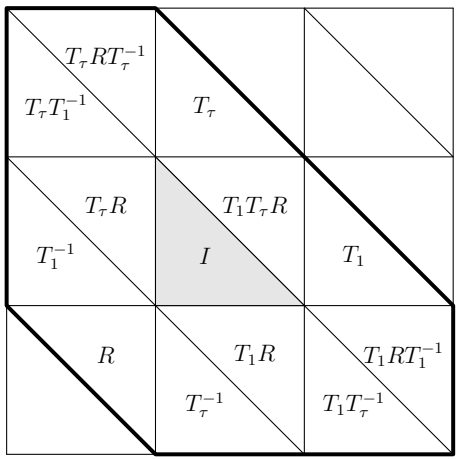}
  \caption{Horizontal tiling in Heisenberg space, using coordinates in the $\R$-basis $1,\tau$}
  \label{fig:tiling}
\end{figure}

We first consider cusp elements. Suppose $\gamma\in\Gamma_\infty$ and
$\gamma(\Omega)\cap\Omega\neq\emptyset$. Then
$\gamma(P)\cap P\neq\emptyset$. Looking at the first component in
Heisenberg space $\C\times\R$, we easily find neighboring triangles in
the tiling, which are obtained by applying elements of the form
$$
\gamma = T_1^j T_\tau^k (T_1T_\tau R)^\epsilon T_v^l,
$$
with $-1,\leq j,k,j+k\leq 1$, $\epsilon=0$ or $1$. A picture of the
corresponding tiling is given in Figure~\ref{fig:tiling}.
For each such element, it is easy to check which values of $l$ give
$\gamma(P)\cap P\neq \emptyset$ (in all cases $-1 \leq l\leq 1$ is
necessary).

Checking the corresponding elements, we find that the only cusp
elements in $T$ (see equation~\ref{eq:T}) are
$$
  R, T_1 R T_1^{-1}, T_\tau R T_\tau^{-1}, T_\tau R, T_1T_\tau R, 
$$
each of order 2.

We will shortly see that $T_\tau R$ and $T_1T_\tau R$ are conjugate in
$\Gamma$, but they are not conjugate to $R$. In particular, we
have~(1) in the following.
\begin{prop}
  Let $\gamma\in\Gamma_\infty$.
  \begin{enumerate}
  \item If $\gamma$ is a non-trivial torsion element, then it is
    conjugate in $\Gamma$ to $R$ or to $T\tau R$.
  \item If $\gamma$ is parabolic but not unipotent, then it is
    conjugate in $\Gamma$ to $T_1RT_v^k$, $RT_v^k$ or $T_\tau RT_v^k$
    for some $k\in\Z$.
  \end{enumerate}  
\end{prop}
Point~(2) follows from the fact that the prism is a fundamental domain
for the action of $\Gamma_\infty$ in Heisenberg space. Indeed,
invariant complex lines for non-unipotent cusp elements can be assumed
to meet the boundary at infinity in one of the four vertical lines
$(0,t)$, $(\frac{1}{2},t)$, $(\frac{\tau}{2},t)$,
$(\frac{1+\tau}{2},t)$ in Heisenberg. The last two are conjugate in
$\Gamma$ (but the first three have distinct conjugacy classes in
$\Gamma_\infty$).

Since our method is not algorithmic for complex reflections, we
list the conjugations needed in order to show that there are indeed
just two conjugacy classes of complex reflections in $\Gamma$ (both of
order 2):
\begin{equation}
  \begin{array}{c}
    (T_1R)^2T_1^{-1} (T_\tau R) ((T_1R)^2T_1^{-1})^{-1} = I\\
    \ [T_1,I] ( T_\tau R ) [T_1,I]^{-1} = T_1 T_\tau R\\
    \left(
    \begin{matrix}
      i\sqrt{7} & 0  & 4\\
      0         & -1 & 0\\
      2         & 0  &-i\sqrt{7}
    \end{matrix}
           \right) = T_1IT_1^{-1} (R) (T_1IT_1)^{-1}.
  \end{array}
  \label{eq:conj-reflections}
\end{equation}
The first two equations show that $T_\tau R$ and $T_1T_\tau R$ are
both conjugate to $I$. The left hand side of the last equation is one
of the complex reflections produced by listing elements of finite
order $\gamma$ with $\Omega\cap\gamma\Omega\neq\emptyset$.

\subsection{Conjugacy classes of torsion elements} \label{sec:other-torsion}

The results below were of course obtained with the help of a computer
(using the methods of section~\ref{sec:algos}). After this paper was
written, we have developed a computer program that performs this
analysis systematically (for general 1-cusped Picard lattices),
see~\cite{picmod}.

For complex reflections, there are precisely two conjugacy classes,
listed in Table~\ref{tab:refl}.
\begin{table}[htbp]
  \centering
  \begin{tabular}{c|c|c}
    Matrix                    & Polar to mirror $v$                 & $\langle v,v\rangle$ \\
    \hline
    $R = \left(\begin{matrix}
      1 & 0  & 0\\
      0 & -1 & 0\\
      0 & 0  & 1
    \end{matrix}\right)$
                              & $(0,1,0)$                                & $1$                  \\
    \hline
    $I = \left(\begin{matrix}
      0 & 0 & 1\\
      0 & -1 & 0\\
      1 & 0 & 0
    \end{matrix}\right)$
                              & $(1,0,1)$                               & $2$                  \\
  \end{tabular}
  \caption{Representatives of conjugacy classes of elements of complex
    reflections (of order 2) in $PU(2,1,\cO_7)$}
  \label{tab:refl}
\end{table}
Note that the two classes can be distinguished by the square norm of
the primitive vector polar to their mirror, which suggests the
following definition.
\begin{dfn}\label{def:jline}
  For $j=1,2$, a $j$-line is a complex line polar to a primitive
  vector $v\in\cO_7$ such that $\langle v,v\rangle = j$.
\end{dfn}

The following result follows from the list of conjugacy classes of
complex reflections in $\Gamma$ (see Table~\ref{tab:refl}).
\begin{prop}
  The group $\Gamma$ acts transitively on the set of primitive vectors
  of square norm $1$ (resp. $2$) in $\cO_K^3$.
\end{prop}
\begin{pf}
  To each $j$-line ($j=1$ or $2$) polar to $v$, the associated complex
  reflection is indeed in $U(J,\cO_7)$, since it is given by
\begin{equation}
  \label{eq:refl}
  R_v(x)=x-2\frac{\langle x,v\rangle}{\langle v,v\rangle} v,
\end{equation}
and $\langle v,v\rangle$ divides 2.  The transitivity now follows from
the fact that there are exactly two conjugacy classes of complex
reflections in $\Gamma$.
\end{pf}

We list the conjugacy classes of torsion elements with isolated fixed
points in Tables~\ref{tab:order2} through~\ref{tab:order7}. For each
class, we give the norm of a primitive vectors representing the fixed
point. We also list the order of the stabilizer of the fixed point
(obtained with the method explained in section~\ref{sec:isotropy}), as
well as the number of 1-lines and 2-lines through that point, see
Definition~\ref{def:jline}.  Recall that being a 1-line (resp. 2-line)
is equivalent to being in the $\Gamma$-orbit of the mirror of $R$
(resp. of $I$).
\begin{table}[htbp]
  \centering
  \begin{tabular}{c|c|c|c|c|c}
    Matrix                    & Fixed pt  $v$                 & $\langle v,v\rangle$  & $|Stab|$    &   1-lines  & 2-lines\\
    \hline
    $(RT_1IT_1^{-1})^2 = \left(\begin{matrix}
      i\sqrt{7} & 0 & 4\\
      0 & 1 & 0\\
      2 & 0 & -i\sqrt{7}
    \end{matrix}\right)$
                             & $(-\bar\tau,0,1)$                 &  $-1$                 &   8       &      2      &      2\\      
    \hline
    $IR = \left(\begin{matrix}
      0 & 0 & 1\\
      0 & 1 & 0\\
      1 & 0 & 0
    \end{matrix}\right)$
                              & $(-1,0,1)$                       &  $-2$                &   4      &      1       &       1\\
    \hline
    $T_1^2I(T_1^{-1}I)^2T_1^2I = \left(\begin{matrix}
      -\bar\tau & \tau & 2\\
      \tau & 2 & \bar\tau\\
      2 & \bar\tau & -\tau
    \end{matrix}\right)$
                            & $(\tau,1,\bar\tau)$               & $-2$                  &  8      &      0       &       4 \\
  \end{tabular}
  \caption{Representatives of conjugacy classes of elements of order 2
    with an isolated fixed point in $PU(2,1,\cO_7)$}
  \label{tab:order2}
\end{table}

\begin{table}[htbp]
  \centering
  \begin{tabular}{c|c|c|c|c|c}
    Matrix                    & Fixed pt $v$                            & $\langle v,v\rangle$ &  $|Stab|$    &   1-lines  & 2-lines\\
    \hline
    $T_1 (I T_1^{-1})^2 (I T_1)^2 = \left(\begin{matrix}
      -1 & \tau & 1\\
      -\bar\tau & 1 & 0\\
      1 & 0 & 0
    \end{matrix}\right)$
                              & $(\tau,1,-\tau)$                           & $-3$                 &   6       &      0        &     3\\
    \hline
    $T_1^2(IT_1^{-1})^2 (IT_1)^2 = \left(\begin{matrix}
        -1 & -\bar\tau & 1\\
        \tau & 1 & 0\\
        1 & 0 & 0
      \end{matrix}\right)$
                              & $(-\bar\tau,1,\bar\tau)$                  & $-3$                 &  6        &    0          &      3\\
    \hline
    $(RT_1)^2 IT_1IT_1^{-1} RT_1IT_1^{-1} = \left(\begin{matrix}
      5 & 0 & -3i\sqrt{7}\\
      0 & -1 & 0\\
      -i\sqrt{7} & 0 & -4
    \end{matrix}\right)$
                            & $\notin\K^3$                   &                       &  6        &    1           &     0\\
  \end{tabular}
  \caption{Representatives of conjugacy classes of elements of order 3
    (isolated fixed point) in $PU(2,1,\cO_7)$}
  \label{tab:order3}
\end{table}

\begin{table}[htbp]
  \centering
  \begin{tabular}{c|c|c|c|c|c}
    Matrix                    & Fixed pt $v$                           & $\langle v,v\rangle$   &   $|Stab|$    &   1-lines  & 2-lines\\
    \hline
    $IT_1^{-1}RT_1 = \left(\begin{matrix}
      0 & 0 & 1\\
      0 & 1 & 2\\
      1 & -2 & -2
    \end{matrix}\right)$
                             & $(-1,-1,1)$                                & $-1$                   &    8       &    2         &    2\\
    \hline
    $(T_1^{-1}I)^2(T_1I)^2 =\left(\begin{matrix}
      0 & 0 & 1\\
      0 & 1 & 1+\bar\tau\\
      1 & -1-\tau & -2
    \end{matrix}\right)$
                           &  $(\bar\tau,-\tau,-\bar\tau)$             & $-2$                   &    8       &    0          &   2+2\\
  \end{tabular}
  \caption{Representatives of conjugacy classes of elements of order 4
    (isolated fixed point) in $PU(2,1,\cO_7)$. The occurrence of $2+2$
    in the last column means that there are four $2$-lines through the
    fixed point, that come in two distinct orbits under the action of the stabilizer.}
  \label{tab:order4}
\end{table}

\begin{table}[htbp]
  \centering
  \begin{tabular}{c|c|c|c|c}
    Matrix                    & Fixed pt $v$                           &  $|Stab|$    &   1-lines  & 2-lines\\
    \hline
    $T_1 R (T_1 I)^2 T_1^{-1} R T_1 I T_1^{-1} = \left(\begin{matrix}
      -5 & 0 & 3i\sqrt{7}\\
      0 & -1 & 0\\
      i\sqrt{7} & 0 & 4
    \end{matrix}\right)$
                             & $\notin \K^3$                   &   6       &        1       &   0\\
  \end{tabular}
  \caption{Representatives of conjugacy classes of elements of order 6
    (isolated fixed point) in $PU(2,1,\cO_7)$}
  \label{tab:order6}
\end{table}

\begin{table}[htbp]
  \centering
  \begin{tabular}{c|c|c|c|c}
    Matrix                    & Fixed pt $v$                           &  $|Stab|$    &   1-lines  & 2-lines\\
    \hline
    $IRT_1=\left(\begin{matrix}
      0 & 0 & 1\\
      0 & 1 & 1\\
      1 & -1 & -\bar\tau
    \end{matrix}\right)$
                              &  $\notin \K^3$                  &     7      &    0          &    0\\
  \end{tabular}
  \caption{Representatives of conjugacy classes of elements of order 7
    (isolated fixed point) in $PU(2,1,\cO_7)$}
  \label{tab:order7}
\end{table}

\section{Study of the mirror stabilizers} \label{sec:mirror-stabs}

\subsection{Mirror of $R$} \label{sec:stabR}

The mirror of $R$ is given by the complex lines in $e_2^\perp$, so it
is quite clear that its stabilizer is isomorphic to
$\{\pm 1\}\times U(1,1,\cO_7)$. In terms of our standard horospherical
coordinates $(z,t,u)$, the mirror is simply given by the points with $z=0$,
i.e. we get a copy of the upper half space $(t,u)\in \R^2$, $u>0$.

The stabilizer can be understood with the same method as
in~\cite{mark-paupert}, only the computations are simpler. Concretely,
to work out the Ford domain for the stabilizer, we restrict to $z=0$
and consider the domain bounded by Cygan spheres for elements
$\gamma\in\Gamma$ only when $\gamma$ preserves $z=0$ (this last is
equivalent to saying that its matrix representative should have
$(1,0,0)$ as an eigenvector, and it implies that $\gamma(\infty)$ have
0 as its second homogeneous coordinate).  Of course there are
infinitely many such elements, but up to the action of
$\Gamma_\infty$, the points $\gamma(\infty)$ must be the ones listed
in~\cite{mark-paupert}.

Also, the cusp elements that preserve $z=0$ consists precisely of
elements in the infinite cyclic group generated by the vertical
translation $T_v$.

It turns out that, up to the action of the group generated by $T_v$,
there is a unique element of depth 1, namely
$$
I = \left(
  \begin{matrix}
    0 & 0 & 1\\
    0 & -1 & 0\\
    1 & 0 & 0
  \end{matrix}\right), 
$$
and a unique element of depth 2, given by
$$
M=\left(
\begin{matrix}
    i\sqrt{7} & 0 & 4\\
    0 & 1 & 0\\
    2 & 0 & -i\sqrt{7}
  \end{matrix}\right).
$$
Moreover, the corresponding Cygan spheres intersect at the (isolated)
fixed point of
$$
M T_v I,
$$
which has order 6 (its cube is actually equal to $R$, so the
transformation acts as an element of order 3 in restriction to the
mirror of $R$).

\begin{figure}[htbp]
  \centering
  \includegraphics[height=8cm]{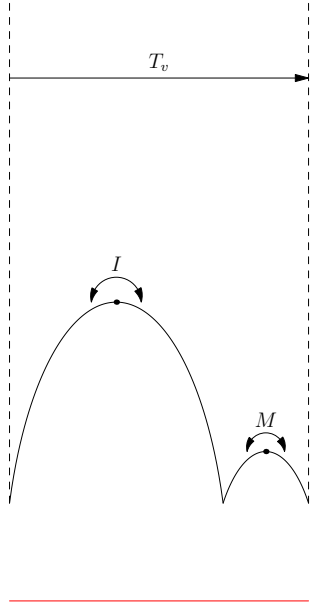}
  \caption{Fundamental domain with side-pairing for the action of the
    stabilizer of the mirror of $R$.}
\label{fig:stabR}
\end{figure}

It follows from the results in~\cite{mark-paupert} that the Cygan
spheres of level 11 or higher do not intersect the Ford domain, so the
determination of the Ford domain for the stabilizer is reduced to
finitely many verifications, and one verifies that the $T_v$-translates
of the two Cygan spheres for $I$ and $M$ actually bound the domain,
and that a fundamental domain for the action of the stabilizer is as
illustrated in Figure~\ref{fig:stabR}.

It is also easy to work out the vertex cycles of the corresponding
Ford polygon, to get a presentation for the fixed point stabilizer of
the form
$$
\langle \iota, \mu, \upsilon | \iota^2, \mu^2, (\mu\upsilon\iota)^3 \rangle.
$$
The stabilizer is a central extension of this group with presentation 
$$
\langle \iota, \mu, \upsilon, \rho | \iota^2, \mu^2, \rho^2, (\mu\upsilon\iota)^3\rho^{-1}, [\rho,\iota],[\rho,\mu],[\rho,\upsilon] \rangle.
$$

The above discussion shows in particular that the mirror stabilizer
has precisely three orbits of points with non-trivial isotropy groups,
represented by
\begin{itemize}
\item The common fixed point of $I$ and $R$ (intersection of one 1-line and one 2-line)
\item The fixed point of $(RT_1IT_1^{-1})^2$ (intersection of two 1-lines and two 2-lines)
\item The fixed point of $MT_vI$ (only one 1-line).
\end{itemize}

\subsection{Mirror of $I$} \label{sec:stabI}

Instead of considering the stabilizer of the mirror of $I$, we will
study the mirror of $T_\tau R$ (which is conjugate to $I$ in
$\Gamma$), because its mirror goes through the point at infinity in
the Siegel half space.

Let $L$ be the mirror of $T_\tau R$, which is given by $v^\perp$ with
$v=(1,-\tau,0)$. The structure of the stabilizer of $L$ in $\Gamma$ is
significantly more difficult to study. Our main result is the
following.
\begin{thm}
  The projection of the stabilizer of $L$ to $PU(1,1)$ is the lattice
  generated by seven elements $r_1,r_2,r_3,r_4$, $s_1,s_2$, with
  presentation
  \begin{equation}\label{eq:pres-stabI}
  \langle r_1,r_2,r_3,r_4,s_1,s_2,t_v\ |\ r_1^2,r_2^3,r_3^2,r_4^2, (s_2^{-1}s_1)^2, s_1^{-1}r_4r_1r_3t_vr_2\ \rangle.
   \end{equation}
   This group has precisely two cusps, corresponding to the cyclic
   groups generated by $t_v$ and $s_2$. The image of the mirror in the
   quotient is a $P^1_\C$ with 2 punctures and 4 orbifold points of
   weight 2.
\end{thm}
Just as in section~\ref{sec:stabR}, $Stab(L)$ is in fact a central
extension of this group, with center generated by the complex
reflection $T_\tau R$ itself (which acts trivially on $L$).

\begin{figure}[htbp]
  \centering
  \includegraphics[height=8cm]{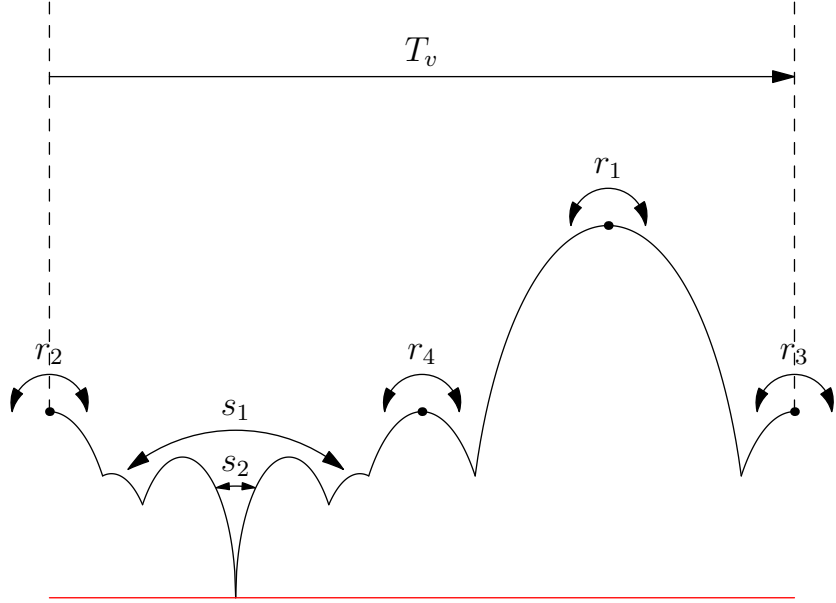}
  \caption{A fundamental domain for the action in $L$ of $Stab(L)$.}\label{fig:ford-stabI}
\end{figure}

We obtained this group by using the computer to list many primitive
vectors $v\in \cO_7^3$ with norm 1 or 2, and keeping only those that
are orthogonal to $v=(1,-\tau,0)$ (so that the corresponding complex
reflections $R_v$ preserves $L$). We then studied the Ford domain (see
Figure~\ref{fig:ford-stabI}) for the group generated by those
reflections, whose side-pairing maps are the above generators.

The group elements $r_1,\dots,r_4$ are complex reflections with mirror
a $2$-line. We describe them by giving vectors $v_1,\dots,v_4$ polar
to their mirror; recall that the matrix of the reflection $R_v$ fixing
$v^\perp$ can be obtained by using formula~\eqref{eq:refl}.
\begin{equation}
  \label{eq:mirrors}
  \begin{array}{c}
    v_1 = (1,1,\bar\tau)\\
    v_2 = (-i\sqrt{7},\tau,2)\\
    v_3 = (i\sqrt{7},\tau,2)\\
    v_4 = (0,1,\bar\tau).
  \end{array}
\end{equation}
For the other two, we give the full matrices, namely
\begin{equation}
  \label{eq:others}
  s_1 =
  \left(
    \begin{matrix}
      -\tau-1 & \tau-2     & \bar\tau + 2\\
      3\tau   &  4         & -5\\
      6       &  3\bar\tau & 5\tau-4
    \end{matrix}\right),\quad
  s_2 =
  \left(
    \begin{matrix}
      \tau-3       & i\sqrt{7}      & -i\sqrt{7}\\
      \bar\tau+3   &   1-i\sqrt{7}  & i\sqrt{7}\\
      -2i\sqrt{7}  &  -\tau-3       & \tau+4
    \end{matrix}\right). 
\end{equation}
It is easy to check that $s_2$ is a parabolic element, whose fixed
point is represented by the primitive vector $(-1,1,\bar\tau)$.

A fundamental domain for the action of the stabilizer of $L$ is
described in Figure~\ref{fig:ford-stabI}. From this, one easily deduces a
presentation of the stabilizer of $L$ (modulo the fixed point
stabilizer of $L$, which is a group of order 2), by using the
Poincar\'e polygon theorem. Indeed, the relations that occur in the
presentation of equation~\eqref{eq:pres-stabI} are the cycle relations
coming from the vertices of the fundamental domain.

A priori the above group is only the (projection to $PU(1,1)$ of the)
subgroup of $Stab(L)$ generated by all the complex reflections in
$Stab(L)$, but in fact this gives full stabilizer, because any
holomorphic symmetry of our domain would have to exchange the two
cusps (but it is easy to see that there is no $A\in U(2,1,\cO_7)$ that
exchanges $(1,0,0)$ and $(-1,1,\bar\tau)$).

\section{A 2-generator presentation} \label{sec:2gen-pres}

\begin{thm}\label{thm:2gen-pres}
  The group $\Gamma_7=PU(2,1,\cO_7)$ has the following presentation
  {\scriptsize
  \begin{equation*}
    \label{eq:2gen-pres}
    \langle\ a,b,c,d\ |\ a^7, b^2, c^6, (ad^2)^4, (c^{-2}d^2)^4, (cd^{-1}c^2d^{-2})^3, (cd^{-2}c^2d^{-1})^3, (d^2c^{-1}a^{-2}d^3c^2a^{-2})^2, c^{-1}ab,d^{-1}ba\ \rangle
  \end{equation*}}
  Moreover, every torsion element in the group is conjugate to an
  element that occurs in one of these power relators.
\end{thm}
Note that we write this as a 4-generator presentation to help
readability, but it should be clear (by looking at the last two
relators) that $a,b$ generate the group, and this can clearly be
turned into a 2-generator presentation (this is useful to speed up
computations using group-theory software like GAP or Magma).

An explicit isomorphism $\phi:G\rightarrow PU(2,1,\cO_7)$ extends
$\phi(a)=A$, $\phi(b)=B$ where
$$
A=\left(
  \begin{matrix}
    -\tau-2 & i\sqrt{7} & i\sqrt{7}\\
    -1      &    1      &   0\\
    \tau-1  &    1      &   1
  \end{matrix}
\right),\quad
B = \left(
  \begin{matrix}
    1 & \bar\tau & -1\\
    0 & -1       & \tau\\
    0 & 0        &   1
  \end{matrix}\right).
$$

The fact that $\phi(a)=A$, $\phi(b)=B$ extends to a group homomorphism
(still denoted by $\phi$) follows from explicit matrix
computation. The fact that this extension is an isomorphism follows
from comparison of $G$ with the presentation for $PU(2,1,\cO_7)$ given
in~\cite{mark-paupert} using the \verb|SearchForIsomorphism| command
in Magma.

We would like to find an element of each torsion conjugacy class
expressed as an explicit word in $a,b$. We start by gathering geometric
data for the obvious elements that occur in the presentation in
equation~\eqref{eq:2gen-pres}, see Table~\ref{tab:obvious-torsion}.
\begin{figure}
\begin{tabular}[htbp]{c|c|c|c|c}
Order &  Group element        & Fixed point $v$            & $\langle v,v\rangle$  &  Other descr.\\
  \hline
2     &  $b$                  & $(1,-\tau,0)$              & $2$                   & $T_\tau R$\\
2     &  $(ba)^3=d^3=a^{-1}c^3a$            & $(\tau,0,1)$                & $1$                   & $T_1IT_1^{-1} R T_1IT_1^{-1}$\\
  \hline
2     &  $((aba)^{-1}babab)^2=(d^{-2}c^2)^2$  & $(\tau,1,\bar\tau)$      & $-2$                  & $\#3$\\
2     &  $(ababa)^2=(ad^2)^2$          & $(\tau+1,1,\bar\tau)$      & $-1$                  & $A_2(\#1)A_2^{-1}$\\                            
  \hline
3     &  $(ba)^2=d^2$             & not in $\Q(i\sqrt{7})$     &                       & $T_1(IT_1^{-1}R)^3$\\
3     &  $[b,a^{-1}babab]=c^{-1}d^2c^{-2}d$    & $(3+i\sqrt{7},1,\bar\tau)$ & $-3$                  & $T_v I(T_\tau J)ITv^{-1}$\\
  \hline
4     &  $(aba)^{-1}babab=d^{-2}c^2$      & $(\tau,1,\bar\tau)$      & $-2$                  & $IT_1^{-1}(IT_1)^2IT_1^{-1}$\\
4     &  $ababa=ad$              & $(\tau+1,1,\bar\tau)$      & $-1$                  & $T_1I(T_1^{-1}I)^2T_1IRT_1IT_1^{-1}$\\                            
  \hline
6     &  $ab=c$                 & not in $\Q(i\sqrt{7})$     &                       & $RT_1IR(T_1I)^2T_1^{-2}$\\
  \hline
7     &  $a$                  & not in $\Q(i\sqrt{7})$     &                       & $T_1RT_1IRT_1I$ 
\end{tabular}
\caption{Obvious elements of finite order, obtained from the
  presentation~\eqref{eq:2gen-pres}.}\label{tab:obvious-torsion}
\end{figure}
In order to get all torsion conjugacy classes (up to replacing any
element to a power that generates the same finite cyclic group), we
need to write a conjugate of $IR$ and a conjugate of
$T_{\bar\tau}^{-1}J$ as explicit words in $a$ and $b$. This can be
done using the above explicit isomorphism between the presentations
(or alternatively by geometric methods using a suitable Dirichlet
domain for the group).

The results are give in Table~\ref{tab:other-torsion}.
\begin{figure}
\begin{tabular}[htbp]{c|c|c|c|c}
  Order   &  Group element                                     & Fixed point $v$            & $\langle v,v\rangle$  &  Other descr.\\
  \hline
  2     & $(aba)^{-1}  ( d^2c^{-1}a^{-2}d^3c^2a^{-2} )  aba$                 &     $(1,0,-1)$ &   $-2$             &    $J=RI$\\
  3     &  $a^{-1}ba^{-1}bababa^{-1}bab=d^{-2}c^2d^{-1}c$                     & $(\tau+1,\bar\tau,-\tau)$  & $-3$                  & $T_1I(T_{\bar\tau}^{-1}J)IT_1^{-1}$\\
\end{tabular}
\caption{Non-obvious elements of finite order.}\label{tab:other-torsion}
\end{figure}

\section{Torsion-free subgroups} \label{sec:torsion-free-subgroup}

We now use the results of the previous sections to exhibit an explicit
torsion-free subgroup of $\Gamma$. Computer code to find this subgroup
(as well as many other torsion-free subgroups of $PU(J,\cO_d)$ for
other values of $d$) is available in~\cite{picmod}.

Recall that $\Gamma$ contains
elements of order $2,3,4,6$ and $7$, so the index of any torsion-free
subgroup has to be a multiple of $\textrm{lcm}(\{4,6,7\})=84$.

Consider the ideal $I=\langle i\sqrt{7}\rangle$ in $\cO_7$, that
satisfies $\cO_7/I\equiv \F_7$, and let
$\phi:\Gamma\rightarrow GL(3,\F_7)$ be the corresponding homomorphism
(reduce all coefficients modulo $I$). We use the matrices $A$ and $B$
from section~\ref{sec:2gen-pres}, and compute
$$
\phi(A)=
\left(
  \begin{matrix}
    1  & 0 & 0\\
    -1 & 1 & 0\\
    3 & 1 & 1
  \end{matrix}
\right),\quad
\phi(B)=
\left(
  \begin{matrix}
    1  & 4 & -1\\
    0 & -1 & 4\\
    0 & 0 & 1
  \end{matrix}
\right).
$$

\begin{prop}
  The group $\Gamma_0=\textrm{Ker}(\phi)$ is a torsion-free subgroup
  of index 336 in $\Gamma$, which is torsion-free at infinity as well.
\end{prop}

\begin{pf}
  It is easy to verify that the subgroup generated by $\phi(A)$ and
  $\phi(B)$ has order 336 (this is most conveniently done with group
  theory software, say GAP or Magma), which gives the statement about
  the index.
  
  The claim about $\Gamma_0$ being torsion-free amounts to verifying
  that every representative $\gamma$ of a conjugacy class of torsion
  elements has the property that $\phi(\gamma)$ has the same order as
  $\gamma$.

  The claim about $\Gamma_0$ being torsion-free at infinity is
  equivalent to checking that $\phi(T_1R)\phi(T_v)^k$ is non-trivial
  for every $k$. In fact $\phi(T_v)$ has order 4, so it it enough to
  check this for $k=0,1,2,3$.
  
  All these verifications are checked by direct computations in
  $GL(3,\F_7)$.  
\end{pf}

\begin{rk}
  \begin{enumerate}
  \item Applying the same construction, but replacing $I$ by the ideal
    $I'=\langle \tau\rangle$ (which satisfies $\cO_7/I'\equiv \F_2$)
    gives a subgroup of index 168, which is not torsion-free.
  \item Using Magma, one can check that $\Gamma$ has exactly two
    normal subgroups of index 336, and only one of them is
    torsion-free.
  \item The group of order 336 $G=\textrm{Im}(\phi)$ is not isomorphic
    to the Shephard-Todd $G_{24}$ that was used
    in~\cite{deraux-klein}. In fact $G$ has trivial center, whereas
    $G_{24}$ has center of order 2.
  \item Still using Magma, one can check that the subgroup $\Gamma_0$
    has finite Abelianization isomorphic to $(\Z/7\Z)^8$. Using
    methods similar to the ones used in the companion computer file
    of~\cite{stover-cusps}, one can also count the number of cusps of
    $\Gamma_0$. It turns out it has 24 cusps, each with
    self-intersection $-7$; more details on this are given
    in~\cite{deraux-xu} (see also the computer code~\cite{picmod}
    related to that paper).
  \end{enumerate}
\end{rk}


\begin{thebibliography}{10}

\bibitem{beardon}
Alan~F. {Beardon}.
\newblock {\em {The geometry of discrete groups}}, volume~91 of {\em {Grad.
  Texts Math.}}
\newblock Springer, New York, NY, 1983.

\bibitem{picmod}
Martin Deraux.
\newblock gitlab project pic-mod.
\newblock \url{https://plmlab.math.cnrs.fr/deraux/pic-mod}.

\bibitem{deraux-family}
Martin {Deraux}.
\newblock {A 1-parameter family of spherical CR uniformizations of the figure
  eight knot complement}.
\newblock {\em {Geom. Topol.}}, 20(6):3571--3621, 2016.

\bibitem{deraux-klein}
Martin {Deraux}.
\newblock {Non-arithmetic lattices and the Klein quartic}.
\newblock {\em {J. Reine Angew. Math.}}, 754:253--279, 2019.

\bibitem{deraux-xu}
Martin Deraux and Mengmeng Xu.
\newblock Torsion in 1-cusped picard modular groups.
\newblock \verb|https://arxiv.org/abs/2205.03037|.

\bibitem{falbel-francsics-parker}
Elisha {Falbel}, G\'abor {Francsics}, and John~R. {Parker}.
\newblock {The geometry of the Gauss-Picard modular group}.
\newblock {\em {Math. Ann.}}, 349(2):459--508, 2011.

\bibitem{falbel-parker}
Elisha {Falbel} and John~R. {Parker}.
\newblock {The geometry of the Eisenstein-Picard modular group}.
\newblock {\em {Duke Math. J.}}, 131(2):249--289, 2006.

\bibitem{ghoshouni-heydarpour-11}
Mahboubeh {Ghoshouni} and Majid {Heydarpour}.
\newblock {A Set of Generators for the Picard Modular Group in the Case d=11}.
\newblock {\em {Iran J. Sci. Technol. Trans. Sci.}}, 44:1469--1475, 2020.

\bibitem{ghoshouni-heydarpour-2}
Mahboubeh {Ghoshouni} and Majid {Heydarpour}.
\newblock {A set of generators for the Picard modular group
  \({SU}(2,1,\mathcal{O}_2)\)}.
\newblock {\em {Proc. Indian Acad. Sci., Math. Sci.}}, 130(26), 2020.

\bibitem{goldman-book}
William~M. {Goldman}.
\newblock {\em {Complex hyperbolic geometry}}.
\newblock Oxford: Clarendon Press, 1999.

\bibitem{humphreys}
James~E. {Humphreys}.
\newblock {\em {Arithmetic groups}}, volume 789 of {\em {Lect. Notes Math.}}
\newblock Springer, 1980.

\bibitem{kim-parker}
Inkang Kim and John~R. Parker.
\newblock Geometry of quaternionic hyperbolic manifolds.
\newblock {\em Math. Proc. Camb. Philos. Soc.}, 135(2):291--320, 2003.

\bibitem{mark-paupert}
Alice Mark and Julien Paupert.
\newblock Presentations for cusped arithmetic hyperbolic lattices.
\newblock \verb|https://arxiv.org/abs/1709.06691|.

\bibitem{parker-will}
John~R. {Parker} and Pierre {Will}.
\newblock {A complex hyperbolic Riley slice}.
\newblock {\em {Geom. Topol.}}, 21(6):3391--3451, 2017.

\bibitem{paupert-will}
Julien {Paupert} and Pierre {Will}.
\newblock {Real reflections, commutators, and cross-ratios in complex
  hyperbolic space}.
\newblock {\em {Groups Geom. Dyn.}}, 11(1):311--352, 2017.

\bibitem{polletta}
David {Polletta}.
\newblock {Presentations for the Euclidean Picard modular groups}.
\newblock {\em {Geom. Dedicata}}, 210:1--26, 2021.

\bibitem{stover-cusps}
Matthew {Stover}.
\newblock {Cusps of Picard modular surfaces}.
\newblock {\em {Geom. Dedicata}}, 157:239--257, 2012.

\bibitem{zhao}
Tiehong {Zhao}.
\newblock {Generators for the Euclidean Picard modular groups}.
\newblock {\em {Trans. Am. Math. Soc.}}, 364(6):3241--3263, 2012.

\bibitem{zink}
Thomas {Zink}.
\newblock {\"uber die Anzahl der Spitzen einiger arithmetischer Untergruppen
  unit\"arer Gruppen}.
\newblock {\em {Math. Nachr.}}, 89:315--320, 1979.

\end{thebibliography}
\end{document}